\documentclass[10pt,twoside]{article}
\usepackage{amsmath,amssymb,amsfonts,amsthm}

\makeatletter

\@addtoreset{equation}{section} \makeatother
\newtheorem{theorem}{Theorem}[section]
\newtheorem{lemma}[theorem]{Lemma}
\newtheorem{proposition}[theorem]{Proposition}

\theoremstyle{definition}
\newtheorem{definition}[theorem]{Definition}

\newtheorem{remark}[theorem]{Remark}

\newcommand{\e}{\mathrm{e}}

\newcommand{\g}{\ensuremath{\langle\cdot,\cdot\rangle}} 

\newcommand{\wk}{\rightharpoonup}
\newcommand{\then}{\Longrightarrow}
\newcommand{\m}{{\mathcal M}}
\newcommand{\J}{{\mathcal J}}
\newcommand{\N}{{\mathbb N}}

\newcommand{\R}{{\mathbb R}}
\newcommand{\acca}{{\mathcal H}}
\newcommand{\elle}{{\mathcal L}}

\begin{document}

\title{\bf\Large A note on geodesic connectedness\\
of G\"odel type spacetimes\thanks{
Par\-tia\-lly sup\-por\-ted by the Spa\-ni\-sh--Ita\-lian
Ac\-ci\'on In\-te\-gra\-da HI2008.0106/A\-zio\-ne In\-te\-gra\-ta
Italia--Spagna IT09L719F1.}}

\author{{\scshape\large R. Bartolo$^\natural$, A.M. Candela\footnote{Partially
supported by Research Funds ``Fondi di Ateneo 2010'' of the Universit\`a degli
Studi di Bari ``Aldo Moro''.} , J.L. Flores\footnote{
Partially supported by the Spanish Grants MTM2010-18099 (MICINN) and
P09-FQM-4496 (J. Andaluc\'{i}a), both with FEDER funds.}}\\
{\small $^\natural$Dipartimento di Matematica, Politecnico di Bari,}\\
{\small Via E. Orabona 4, 70125 Bari, Italy}\\
{\small \texttt{r.bartolo@poliba.it}}\\
{\small $^\dagger$Dipartimento di Matematica, Universit\`a degli Studi di Bari ``A. Moro'',}\\
{\small Via E. Orabona 4, 70125 Bari, Italy}\\
{\small \texttt{candela@dm.uniba.it}}\\
{\small $^\ddagger$Departamento de \'Algebra, Geometr\'{\i}a y Topolog\'{\i}a,}\\
{\small Facultad de Ciencias, Universidad de M\'alaga,}\\
{\small Campus Teatinos, 29071 M\'alaga, Spain}\\
{\small \texttt{floresj@agt.cie.uma.es}}}

\date{}

\maketitle

\begin{abstract}
In this note we reduce the problem of geodesic connectedness in a wide
class of G\"odel type spacetimes to the search of critical points of a functional
naturally involved in the study of geodesics in standard static
spacetimes. Then, by using some known accurate results on the
latter, we improve previous results on the former.
\end{abstract}


\noindent
{\bf Corresponding author:}

\noindent

Rossella Bartolo

\noindent

Dipartimento di Matematica, Politecnico di Bari

\noindent

Via E. Orabona 4, 70125 Bari, Italy

\noindent

e-mail: r.bartolo@poliba.it; fax: +39 080 5963612

\newpage

\section{Introduction}
In the last years the study of geodesic connectedness in some
classes of Lorentzian manifolds has been carried out
systematically by using variational methods (see \cite{CS2} and
references therein). This is the case of static and wave type
spacetimes, where the variational methods
yield optimal results by reducing the
strongly indefinite action functional
to a subtler Riemannian one (cf. \cite{bcfs, CFS2}).

In \cite{CS} a similar approach is applied to G\"odel type
spacetimes (here, Definition \ref{type}), but without providing optimal
results (here, Theorem \ref{connect}). In fact, in this case the
variational result does not cover the classical G\"odel Universe,
whose geodesic connectedness is proved instead by a direct
integration of the corresponding geodesic equations (see
\cite[Section 4]{CS}).

The aim of this paper is to improve meaningfully the result in
\cite{CS} by a more careful application of variational methods.
Indeed, we can deal with a functional similar to that one defined
for standard static spacetimes (Section \ref{static}). Then, by
using the accurate estimates in \cite{bcfs}, we provide a
substantial weakening of the boundedness assumptions about the
metric coefficients in \cite{CS} (see Theorems \ref{connect1} and \ref{connect2}).
Unfortunately, the fact that our theorems do not cover the
classical G\"odel Universe seems to indicate that a sharp result
cannot be reached only by using variational methods.

The paper is organized as follows: in Section \ref{tools}, for the
reader's convenience, we recall some definitions about the
variational setting in the problem of geodesic connectedness; in
Sections \ref{goedel} and \ref{static} we introduce the notions,
and related results, of G\"odel type and static spacetime,
respectively; finally, in Section \ref{s2} we prove our main
theorems.

\section{Variational setting}\label{tools}
In order to state our main results, firstly we recall some notations useful
 for the variational setting.

Taking a connected, finite--dimensional semi--Riemannian
manifold $(\m,g)$, let $H^1(I,\m)$ be the set of curves $z : I \to \m$, $I =
[0,1]$, such that for any local chart $(U,\varphi)$ of $\m$,
with $U \cap z(I) \ne \emptyset$, the curve $\varphi \circ z$
belongs to the Sobolev space $H^1(z^{-1}(U),\R^n)$, $n = {\rm
dim}\, \m$. Then, $H^1(I,\m)$ is equipped with a structure of
infinite--dimensional manifold mo\-del\-led on the Hilbert space
$H^1(I,\R^n)$. For any $z \in H^1(I,\m)$ the tangent space of
$H^1(I,\m)$ at $z$ can be written as follows:
\[
T_zH^1(I,\m) = \{\zeta \in H^1(I,T\m) : \zeta(s) \in
T_{z(s)}\m \; \mbox{for all $s \in I$}\},
\]
where $T\m$ is the tangent bundle of $\m$.

If $\m$ splits globally in the product of two semi--Riemannian
mani\-folds $\m_1$ and $\m_2$, i.e. $\m=\m_1\times\m_2$, then
\[
H^1(I,\m)\equiv H^1(I,\m_1)\times H^1(I,\m_2)
\]
and $T_zH^1(I,\m)\equiv T_{z_{1}}H^1(I,\m_1)\times
T_{z_{2}}H^1(I,\m_2)$ for all $z=(z_1,z_2)\in\m$.

On the other hand, if $(\m_0,\langle\cdot,\cdot\rangle_R)$ is a
complete Riemannian manifold, it can be smoothly and isometrically
embedded in an Euclidean space $\R^N$ (cf. \cite{mu}). Hence,
$H^1(I,\m_0)$ is a submanifold of the Hilbert space $H^1(I,\R^N)$.
In this case, we denote by $d(\cdot,\cdot)$ the distance
induced on $\m_0$ by its Riemannian metric
$\langle\cdot,\cdot\rangle_{R}$, i.e.
\[
d(x_p,x_q)\ :=\ \inf \left\{  \int_a^b \sqrt{\langle\dot x,\dot
x\rangle_R} \; {\rm d}s : \; x \in A_{x_p,x_q}\right\},
\]
where $x \in A_{x_p,x_q}$ if $x: [a,b]\rightarrow\m_0$ is any
piecewise smooth curve in $\m_0$ joining $x_p, x_q\in\m_0$.

Given $z_p$, $z_q \in \m$, let us consider
\[
\Omega^1(z_p,z_q) = \{z\in H^1(I,\m) : z(0) = z_p ,\; z(1) = z_q\},
\]
which is a (complete if $\m$ is complete) submanifold of
$H^1(I,\m)$, having tangent space at any $z \in
\Omega^1(z_p,z_q)$ described as
\[
T_z\Omega^1(z_p,z_q) = \{\zeta \in T_zH^1(I,\m) : \zeta(0) = 0 =
\zeta(1)\}.
\]
Moreover, for any $l_p$, $l_q \in \R$, let us denote
\[
W(l_p,l_q) = \{l \in H^1(I,\R) : l(0) = l_p\ ,\; l(1) = l_q\} .
\]
Clearly,
\[
W(l_p,l_q) = H^1_0(I,\R) + \bar l ,
\]
with $H^1_0(I,\R) = \{l \in H^1(I,\R) : l(0) = 0 = l(1)\}$, $\bar
l : s \in I \mapsto (1 - s) l_p + s l_q \in \R$. Hence,
$W(l_p,l_q)$ is a closed affine submanifold of the Hilbert space
$H^1(I,\R)$ with tangent space
\[
T_lW(l_p,l_q) = H^1_0(I,\R) \quad \hbox{for every $l \in W(l_p,l_q)$.}
\]
At last, let us recall a classical variational princip
le:
if  $(\m,g)$ is a
semi--Riemannian manifold, then $\bar z:I\rightarrow \m$
is a geodesic joining two points $z_p,z_q\in\m$ if and only if
$\bar z\in \Omega^1(z_p,z_q)$ is a critical point of the action functional
\begin{equation}\label{act}
f(z) = \frac 12\ \int_0^1 g(z)[\dot z,\dot z]\; {\rm d}s \qquad
\hbox{on $ \Omega^1(z_p,z_q)$.}
\end{equation}

\section{G\"odel type spacetimes and statement of the main theorems}\label{goedel}

The classical {\em G\"odel Universe} is an exact solution of the
Einstein's field equations in which the matter takes the form of a
rotating pressure--free perfect fluid. Matematically, it is
modelled by $\R^4$ equipped with metric
\[
d s^2 = d x_1^2 + d x_2^2 - \frac 12\ \e^{2 \sqrt 2 \omega x_1}
dy
^2 - 2\ {\e}^{\sqrt 2 \omega x_1} dy dt - dt^2,\quad x =
(x_1,x_2) \in \R^2,
\]
where $\omega > 0$ is the magnitude of the vorticity of the flow
(see \cite{HE}). It can be proved that $(\R^{4},ds^{2})$ is a
geodesically connected Lorentzian manifold (cf. \cite[Section 4]{CS}).

A natural generalization of the classical G\"odel Universe is
introduced in \cite{CS} as follows.
\begin{definition}\label{type}

{\rm A Lorentzian manifold $(\m, \langle\cdot,\cdot\rangle_L)$ is
a {\em G\"odel type spacetime} if there exists a smooth
(connected) finite--dimensional Riemannian manifold
$(\m_0,\langle\cdot,\cdot\rangle_R)$ such that $\m = \m_0 \times
\R^2$ and the metric $\langle\cdot,\cdot\rangle_L$ is described
as:
\begin{equation}
\label{metric} \langle\cdot,\cdot\rangle_L\ =\
\langle\cdot,\cdot\rangle_R + A(x) dy^2 + 2 B(x) dy dt  - C(x)
dt^2 ,
\end{equation}
where $x \in \m_0$, the variables $(y,t)$ are the natural
coordinates of $\R^2$ and $A$, $B$, $C$ are $C^1$ scalar fields on
$\m_0$ satisfying
\begin{equation}
\label{stima1} \acca(x)\ =\ B^2(x) + A(x) C(x) > 0\quad \mbox{for
all $x \in \m_0$.}
\end{equation}}
\end{definition}

Note that condition \eqref{stima1} implies that metric
\eqref{metric} is Lorentzian. Furthermore, to this class
it belongs not only the G\"odel Universe, just taking
\[
\begin{split}
&(\m_{0},\langle\cdot,\cdot\rangle_R) = (\R^{2}, d x_1^2 + d x_2^2),\\
&A(x) = - {\e}^{2 \sqrt 2 \omega x_1}/2,\quad B(x) = - {\e}^{\sqrt
2 \omega x_1},\quad C(x) \equiv 1,
\end{split}
\]
but also other physically relevant models of Lorentzian manifolds,
as some stationary spacetimes, the Kerr--Schild spacetime or some
warped product spacetimes (see \cite{CS} and references therein).

\begin{remark}
{\rm If the product $A(x) C(x)$ is strictly positive on $\m_0$,
the corresponding G\"odel type spacetime reduces to a standard
stationary one, and its geodesic connectedness has been deeply
studied in previous works (see, for example, \cite{bcf,CFS08}).}
\end{remark}

Every G\"odel type spacetime $(\m, \langle\cdot,\cdot\rangle_L)$
admits two Killing vector fields, $\partial_y$ and $\partial_t$,
which are not necessarily timelike. At a first glance, the search
of geodesics for these spacetimes can be handled in the same
manner as in the static case: in fact, geodesics are critical
points of the corresponding action functional which, as in the
static case, becomes equivalent to a suitable simpler
``Riemannian'' one (cf. \cite[Proposition 2.2]{CS}; here Section
\ref{s2}).

Taking a G\"odel type spacetime $(\m=\m_0\times\R^2, \langle\cdot,\cdot\rangle_L)$
according to Definition \ref{type}, for each $x\in H^1(I,\m_0)$
let us introduce the following notations:
\begin{eqnarray}
& &a(x) = \int_0^1 \frac{A(x)}{\acca(x)} \; {\rm d}s ,\;\; b(x) =
\int_0^1 \frac{B(x)}{\acca(x)}\; {\rm d}s , \;\; c(x) = \int_0^1
\frac {C(x)}{\acca(x)}\; {\rm d}s,
\label{definition1}\\
& &\elle(x) = b^2(x) + a(x) c(x).\nonumber
\end{eqnarray}
Then, the following result on geodesic connectedness in
$(\m,\langle\cdot,\cdot\rangle_L)$ is obtained (see \cite[Theorem
1.3]{CS}, \cite[Theorem 1.3]{CS1}):
\begin{theorem}
\label{connect} Let $(\m= \m_0 \times
\R^2,\langle\cdot,\cdot\rangle_L)$ be a G\"odel type spacetime
such that
\begin{itemize}
\item[$(h_1)$] $(\m_0,\langle\cdot,\cdot\rangle_R)$ is a complete
Riemannian manifold; \item[$(h_2)$] $|\elle(x)| > 0$ for all $x\in
H^1(I,\m_0)$; \item[$(h_3)$] there exist $k_1$, $k_2$, $k_3 > 0$
such that
\[
\left|\frac{a(x)} {\elle(x)}\right| \le k_1 ,\quad
\left|\frac {b(x)} {\elle(x)}\right| \le k_2 ,\quad
\left|\frac {c(x)} {\elle(x)}\right| \le k_3 \quad
\mbox{for all $x \in H^1(I,\m_0)$.}
\]
\end{itemize}
Then, $(\m,\langle\cdot,\cdot\rangle_L)$ is geodesically
connected.
\end{theorem}

Even if the hypotheses of this theorem are not optimal, the
counterexample in \cite[Appendix B]{CS} shows that they are
reasonable. However, the boundedness assumptions in
Theorem \ref{connect} can be improved considerably. As a matter of fact,
preserving $(h_1)$, assumption $(h_3)$ can be replaced by heavily weaker
hypotheses which directly involve the coefficients of the
Lorentzian metric on $\m$. In fact, the main aim of this paper is
to prove the following results:

\begin{theorem}
\label{connect1} Let $(\m= \m_0 \times
\R^2,\langle\cdot,\cdot\rangle_L)$ be a G\"odel type spacetime
such that $(h_1)$ holds. Moreover, assume that
\begin{itemize}
\item[$(h'_2)$] there exists $\nu > 0$ such that $\elle(x) \ge \nu
> 0$ for all $x \in H^1(I,\m_0)$;
\item[$(h'_3)$] $A(x)-C(x)>0$ for all $x \in \m_0$,
and there exist $\lambda \ge 0$, $k \in \R$ and a point $\bar x\in
\m_0$ such that the (positive) map
\[
\gamma: x\in \m_0\ \mapsto\ \frac{\acca(x)}{A(x)-C(x)} \in \R
\]
satisfies
\begin{equation}\label{ipotesi1}
\gamma(x) \le \lambda d^2(x,\bar x) + k \quad \mbox{for all $x \in
\m_0$.}
\end{equation}
\end{itemize}
Then, $(\m,\langle\cdot,\cdot\rangle_L)$ is geodesically
connected.
\end{theorem}

\begin{theorem}
\label{connect2} Let $(\m= \m_0 \times
\R^2,\langle\cdot,\cdot\rangle_L)$ be a G\"odel type spacetime
such that $(h_1)$ holds. Moreover, assume that
\begin{itemize}
\item[$(h''_2)$] there exists $\nu > 0$ such that $\elle(x) \le
-\nu < 0$ for all $x \in H^1(I,\m_0)$; \item[$(h''_3)$]
$A(x)-C(x)<0$ for all $x \in \m_0$.
\end{itemize}
Then, $(\m,\langle\cdot,\cdot\rangle_L)$ is geodesically
connected.
\end{theorem}

\section{Static spacetimes}\label{static}

The proof of Theorem \ref{connect1} relies on some results of variational
nature coming from the study of geodesics in {\sl static
spacetimes}, i.e. Lorentzian mani\-folds endowed with an
irrotational timelike Killing vector field. This section is
dedicated to recall these statements.

\begin{definition}{\rm A Lorentzian manifold
$(\m ,\g_L)$ is a {\em standard static spacetime} if there exists
a smooth (connected) finite--dimensional Riemannian mani\-fold
$(\m_0,\langle \cdot,\cdot\rangle_R)$ such that $\m = \m_0 \times
\R$ and the metric $\langle \cdot,\cdot\rangle_L$ is described as:
\begin{equation}
\label{metrica} \langle\cdot,\cdot\rangle_L = \langle
\cdot,\cdot\rangle_R - \beta(x)\ d t^2,
\end{equation}
with $x \in \m_0$, $t$ the natural coordinate of $\R$ and $\beta$
a smooth strictly positive scalar field on $\m_{0}$.}
\end{definition}
The following two statements about geodesic connectedness in
static spacetimes are well--known:
\begin{itemize}
\item The problem of geodesic connectedness in a (connected)
static spacetime can be reduced to the same problem in a suitable
standard static spacetime (see \cite[Section 2]{bcfs}); \item Two
points $z_{p}=(x_p,t_p)$, $z_{q}=(x_q,t_q)$ of a standard static
spacetime $(\m_0\times\R,\g_L)$ are connected by a geodesic
$z=(x,t)$, which is a critical point of the strongly indefinite
action functional $f$ in \eqref{act}, with $g=\langle\cdot,\cdot
\rangle_L$ as in \eqref{metrica} and
$\Omega^{1}(z_{p},z_{q})=\Omega^1(x_p,x_q)\times W(t_p,t_q)$, if
and only if the functional
\begin{equation}
\label{newfunc} J(x) =\frac 12\int_0^1 \langle \dot x,\dot
x\rangle_R {\rm d}s-\frac{\Delta_t^2}{2}\left(\int_0^1 \frac
1{\beta(x)}\ {\rm d}s\right)^{-1},
\end{equation}
with $\Delta_t:=t_p-t_q$, admits a critical point on
$\Omega^1(x_p,x_q)$ (see \cite{BFG}). This variational principle
is a consequence of the existence of the Killing vector field
$\partial_t$ on $(\m,\langle\cdot,\cdot\rangle_{L})$, which
implies the constancy of $\langle
\partial_t,\dot z\rangle_L$ along each geodesic $z$ on $\m$.
\end{itemize}

The existence of critical points for functional $J$ in
$\Omega^1(x_p,x_q)$, and thus the geodesic connectedness of
standard static spacetimes, is ensured under different conditions
for the growth of the metric coefficient $\beta$: when $\beta$ is
bounded (cf. \cite{BFG}), when it is subquadratic (e.g., cf.
\cite{CMP}), and when it grows at most quadratically with respect
to the distance $d(\cdot,\cdot)$ induced on $\m_0$ by its
Riemannian metric $\langle\cdot,\cdot\rangle_{R}$. More precisely
(\cite[Theorem 1.1]{bcfs}):

\begin{theorem}\label{maintheorem}
Let $(\m=\m_{0}\times\R,\langle\cdot,\cdot\rangle_L)$ be a
standard static spacetime such that
\begin{itemize}
\item[$(H_1)$] $(\m_0,\langle\cdot,\cdot\rangle_R)$ is a complete
Riemannian manifold, \item[$(H_2)$] the positive function $\beta$
grows at most quadratically at infinity, i.e. there exist $\lambda
\ge 0$, $k \in \R$ and a point $\bar x\in \m_0$ such that
\[
\beta(x) \le \lambda d^2(x,\bar x) + k
\quad \mbox{for all $x \in \m_0$.}
\]
\end{itemize}
Then, $(\m,\langle\cdot,\cdot\rangle_L)$ is geodesically
connected.
\end{theorem}

In this theorem the growth assumption for $\beta$ is optimal. More
precisely, there exists a family of geodesically disconnected
static spacetimes with superquadratic, but arbitrarily close to
quadratic, coefficients $\beta$ (\cite[Section 7]{bcfs}).

In order to prove Theorem \ref{maintheorem}, the following
classical critical point result is applied (e.g., cf.
\cite[Theorem 2.7]{rab}):

\begin{theorem}
\label{min} Assume that $\Omega$ is a complete Riemannian manifold
and ${\mathcal F}$ is a $C^1$ functional on $\Omega$ which
satisfies the Palais--Smale condition, i.e. any sequence $(x_k)_k
\subset \Omega$ such that
\[
({\mathcal F}(x_k))_k\; \mbox{is bounded}\quad\hbox{and}\quad \lim_{k \to
+\infty}{\mathcal F}'(x_k) = 0
\]
converges in $\Omega$, up to subsequences. Then, if ${\mathcal F}$ is bounded
from below, it attains its infimum.
\end{theorem}

In fact, in our case, $(H_1)$ implies that
$\Omega=\Omega^{1}(x_p,x_q)$ is complete for each
$x_p,x_q\in\m_0$. Moreover, the boundedness and the Palais-Smale
conditions for ${\mathcal F}=J$ are ensured by the following
technical result (cf. \cite[Propositions 4.1, 4.3]{bcfs}):
\begin{proposition}
\label{inf} Under the hypotheses of Theorem \ref{maintheorem}, for
each $x_p,x_q\in\m_0$ the functional $J$ on $\Omega^1(x_p,x_q)$ is
\smallskip

\noindent
$\bullet\; $ bounded from below;\\
$\bullet\; $ coercive, i.e.
$\; \displaystyle J(x) \to +\infty \quad \mbox{as}\quad
\|\dot x\|^2:= \int_0^1 \langle\dot x,\dot x\rangle_R\ {\rm d}s \to +\infty$.
\end{proposition}

\section{On functional $\J$ and proofs of the main theorems}

\label{s2}

In \cite{CS} the authors develop a variational principle which
allows one to study the geodesic connectedness of G\"odel type
spacetimes by finding critical points of a suitable functional
$\J$ (see \eqref{functional} below). After recalling this
principle, in the present section we rewrite functional $\J$ to
find a connection with the static functional $J$ in
\eqref{newfunc}. As a consequence, the geodesic connectedness of
certain G\"odel type spacetimes is deduced as a corollary of
Theorem \ref{maintheorem}.

Throughout this section, $\m = \m_0 \times \R^2$ is a G\"odel type
spacetime according to Definition \ref{type}. Fixing $z_p =
(x_p,y_p,t_p)$, $z_q = (x_q,y_q,t_q) \in \m$, with $x_p$, $x_q \in
\m_0$ and $(y_p,t_p)$, $(y_q,t_q) \in \R^2$, from the product
structure of $\m$ and the remarks in Section \ref{tools}, we have
that $\bar z:I\rightarrow \m$ is a geodesic joining $z_p$ to $z_q$
in $\m$ if and only if $\bar z$ is a critical point of the $C^1$
action functional in \eqref{act} with $g=\langle\cdot,\cdot
\rangle_L$ as in \eqref{metric} and $\Omega^1(z_p,z_q)=
\Omega^1(x_p,x_q) \times W(y_p,y_q) \times W(t_p,t_q)$.

One can take advantage of the Killing vector fields $\partial_y$,
$\partial_t$ on $\m$ for proving a new variational principle free from
the strongly indefinite character of $f$ in $\Omega^1(z_p,z_q)$. In fact, for all $s\in
I$ and every $x \in \Omega^1(x_p,x_q)$ such that $\elle(x) \ne 0$,
consider
\begin{eqnarray}
&&\begin{aligned} \phi_y(x)(s)\ :=\ & y_p + \ \frac{\Delta_y\ b(x)
- \Delta_t\ c(x)} {\elle(x)}\
\int_0^s \frac {B(x)} {\acca(x)}\; {\rm d}\sigma\\
&+\ \frac{\Delta_y\ a(x) + \Delta_t\ b(x)} {\elle(x)}\
\int_0^s \frac{C(x)} {\acca(x)}\; {\rm d}\sigma ,
\end{aligned}\label{uai}\\
&&\begin{aligned} \phi_t(x)(s)\ :=\ & t_p - \ \frac{\Delta_y\ b(x)
- \Delta_t\ c(x)} {\elle(x)}\
\int_0^s \frac{A(x)} {\acca(x)}\; {\rm d}\sigma\\
&+\ \frac{\Delta_y\ a(x) + \Delta_t\ b(x)} {\elle(x)}\ \int_0^s
\frac{B(x)}{\acca(x)}\; {\rm d}\sigma
\end{aligned}\label{ti}
\end{eqnarray}
with
\[
\Delta_y := y_q-y_p,\qquad \Delta_t := t_q-t_p.
\]
Standard arguments imply that
\[
\phi_y : \Omega^1(x_p,x_q) \to W(y_p,y_q)\qquad \mbox{and}\qquad
\phi_t : \Omega^1(x_p,x_q) \to W(t_p,t_q)
\]
are $C^{1}$ functions.

Then, one can establish the following proposition (see
\cite[Proposition 2.2]{CS} for further details):
\begin{proposition}
If $x_p$, $x_q \in \m_0$ are such that $|\elle(x)| > 0$ for all $x
\in \Omega^1(x_p,x_q)$,
 then the following statements are equivalent:
\begin{itemize}
\item[{\sl (i)}] $\bar z \in Z$ is a critical point of the action
functional $f$ in \eqref{act}; \item[{\sl (ii)}] setting $\bar z =
(\bar x,\bar y,\bar t)$, we have that $\bar x \in
\Omega^1(x_p,x_q)$ is a critical point of the $C^1$ functional
\begin{equation} \label{functional} \J(x) = \frac 12\ \int_0^1
\langle\dot x,\dot x\rangle_R \; {\rm d}s \ +\ \frac{\Delta_y^2
a(x) + 2 \Delta_y \Delta_t b(x) - \Delta_t^2 c(x)} {2
\elle(x)}
\end{equation}
on $\Omega^1(x_p,x_q)$, and the other components satisfy $\bar y = \phi_y(\bar x)$, $\bar
t = \phi_t(\bar x)$, with $\phi_y$, $\phi_t$ as above.
\end{itemize}
Furthermore,
\begin{equation}\label{newfunct}
\J(x)\ =\ f(x,\phi_y(x),\phi_t(x)) \quad \hbox{for all $x \in
\Omega^1(x_p,x_q)$.}
\end{equation}
\end{proposition}
Now, we are ready to develop the key point of our approach by
writing functional \eqref{functional} in a smarter way.

Giving $x \in \Omega^1(x_p,x_q)$ such that
$|\elle(x)| > 0$, the numerator of the fraction in
\eqref{functional} is a quadratic form that can be rewritten as follows:
\[
\Delta_y^2 a(x) + 2 \Delta_y \Delta_t b(x) - \Delta_t^2 c(x)\ =\
\left( \begin{array}{cc}
\Delta_y & \Delta_t  \\
\end{array} \right)
\ \left( \begin{array}{cc}
a(x) & b(x) \\
b(x) & - c(x)
\end{array} \right)\
\left( \begin{array}{cc}
\Delta_y  \\
\Delta_t
\end{array} \right).
\]
Note that the symmetric matrix
\[
S(x)\ =\
\left( \begin{array}{cc}
a(x) & b(x) \\
b(x) & - c(x)
\end{array} \right),\qquad\hbox{with}\;\;\;
\det S(x) = - \elle(x) \ne 0,
\]
admits two real (non--null) eigenvalues
\begin{eqnarray}
\nonumber
\lambda_\pm(x)\ &= &\ \frac{a(x)-c(x) \pm \sqrt{(a(x)-c(x))^2 + 4\elle(x)}}{2}\\
&=&\ \frac{a(x)-c(x) \pm \sqrt{(a(x)+c(x))^2 +
4b^2(x)}}{2}\label{*}
\end{eqnarray}
which are the solutions of the characteristic equation
\begin{equation}\label{car}
\lambda^2 - (a(x)-c(x)) \lambda - \elle(x)\ =\ 0.
\end{equation}
Moreover, the following relations hold:
\[
\lambda_-(x) \le \lambda_+(x),\quad \lambda_-(x) \lambda_+(x) = -
\elle(x) \ne 0,\quad \lambda_+(x) + \lambda_-(x) = a(x)-c(x).
\]
These eigenvalues are associated to the normalized eigenvectors:
\[
\tilde v_\pm(x)=\frac{v_\pm(x)}{|v_\pm(x)|}\quad \hbox{with}\quad
v_\pm(x)\ =\ \left(\frac{\lambda_\pm(x) + c(x)}{b(x)},1\right)
\quad \hbox{if $b(x)\neq 0$,}
\]
$\tilde v_+(x) = (1,0)$, $\tilde v_-(x) = (0,1)$
if $b(x) = 0$ and $a(x) > -c(x)$ (being $\lambda_+(x) = a(x)$,
$\lambda_-(x) = -c(x)$), or
$\tilde v_+(x) = (0,1)$, $\tilde v_-(x) = (1,0)$
if $b(x) = 0$ and $a(x) < -c(x)$ (being $\lambda_+(x) = -c(x)$,
$\lambda_-(x) = a(x)$).\\
As a consequence, if $D(x)$ is the matrix whose columns are $\tilde v_\pm(x)$, then
\[
D(x)^{-1}=D(x)^T\qquad \hbox{and}\qquad
D(x)^T S(x) D(x)\ =\ \left( \begin{array}{cc}
\lambda_+(x) & 0 \\
0 & \lambda_-(x)
\end{array} \right).
\]
At any case, it results
\begin{eqnarray*}
&&\left( \begin{array}{cc} \Delta_y & \Delta_t
\end{array} \right) S(x) \left( \begin{array}{cc}
\Delta_y  \\
\Delta_t
\end{array} \right)\\
&&\qquad =\ \left( \begin{array}{cc} \Delta_y & \Delta_t
\end{array} \right)D(x) \left( \begin{array}{cc}

\lambda_+(x) & 0 \\
0 & \lambda_-(x)
\end{array} \right)D(x)^T\left( \begin{array}{cc}\Delta_y  \\
\Delta_t
\end{array} \right)\\
&&\qquad =\ \left( \begin{array}{cc} \Delta_+(x) & \Delta_-(x)
\end{array} \right)\left( \begin{array}{cc}
\lambda_+(x) & 0 \\
0 & \lambda_-(x)
\end{array} \right)\left( \begin{array}{cc}
\Delta_+(x)  \\
\Delta_-(x)
\end{array} \right)\\
&&\qquad =\ \lambda_+(x) \Delta_+^2(x) + \lambda_-(x)
\Delta_-^2(x),
\end{eqnarray*}
where $\Delta_+(x):= (\Delta_y\;\Delta_t)\cdot \tilde v_+(x)$ and
$\Delta_-(x):=(\Delta_y\;\Delta_t)\cdot \tilde v_-(x)$. By
definition, we have
\begin{equation}\label{differe}
|\Delta_+(x)|\leq\ \sqrt{\Delta_y^2+\Delta_t^2},\qquad
|\Delta_-(x)|\leq\ \sqrt{\Delta_y^2+\Delta_t^2}.
\end{equation}
Thus, we obtain
\begin{eqnarray}\nonumber
\J(x)\ &=&\ \frac 12\ \|\dot x\|^2
\ +\
\frac{
\left( \begin{array}{cc}
\Delta_y & \Delta_t  \\
\end{array} \right) S(x) \left( \begin{array}{cc}
\Delta_y  \\
\Delta_t
\end{array} \right)}{2\elle (x)}\\
\nonumber
&=&\
\frac 12\ \|\dot x\|^2
\ -\frac{\lambda_+(x) \Delta_+^2(x) + \lambda_-(x)
\Delta_-^2(x)}{2\lambda_+(x) \lambda_-(x)}\\
\label{newact} &=&\ \frac 12\ \|\dot x\|^2 \ -\ \frac 12\
\frac{\Delta_+^2(x)}{\lambda_-(x)}\ -\ \frac 12\
\frac{\Delta_-^2(x)}{\lambda_+(x)}\ .
\end{eqnarray}

In order to discuss the boundedness and growth behavior of $\J$ in
(\ref{newact}), let us focus on equation \eqref{car}. From
Descartes' rule of sign, the following cases may occur:

\vspace{1cm}

\begin{tabular}{||c|c|c|c|c|c||}\hline

&$\elle(x)$ & $a(x)-c(x)$ &  &  &  \\ \hline

$(i)$ & $\elle(x)>0$& $a(x)-c(x)>0$ & $\then$ & $\lambda_-(x) <0$& $\lambda_+(x) >0$ \\ \hline

$(ii)$ & $\elle(x)>0$& $a(x)-c(x)<0$& $\then$ & $\lambda_-(x) <0$&$\lambda_+(x) >0$ \\ \hline

$(iii)$ & $\elle(x)<0$ & $a(x)-c(x)>0$ & $\then$ & $\lambda_-(x) >0$ & $\lambda_+(x) >0$ \\ \hline

$(iv)$ &$\elle(x)<0$& $a(x)-c(x)<0$ & $\then$ &  $\lambda_-(x) <0$&  $\lambda_+(x) <0$ \\ \hline

\end{tabular}

\vspace{1cm}

\begin{remark}
{\rm From (\ref{*}), the equality $\lambda_-(x) = \lambda_+(x)$
occurs when $b(x) = 0$ and $a(x) = -c(x)$, and it implies
$\elle(x)=-c(x)^{2}<0$\footnote{Note that this situation cannot
occur when $B(x) \equiv 0$ on $\m_0$. In fact, under this
condition hypothesis \eqref{stima1} forces $A(x)$, $C(x)$ to have
the same sign, and the same must happen for $a(x)$, $c(x)$, in
contradiction with the equality $a(x) = -c(x)$.}. So, from
previous table, condition $\elle(x)
> 0$ implies $\lambda_-(x) < 0 < \lambda_+(x)$.}
\end{remark}

\begin{lemma}\label{ps1}
Assume that hypothesis $(h_1)$ holds. Fixing $x_p$, $x_q\in\m_0$,
suppose that $\J$ is coercive on $\Omega^1(x_p,x_q)$ and $\nu > 0$
exists such that
\begin{equation}\label{ps3}
|\elle(x)| \ \ge \nu\quad \hbox{for all $x \in
\Omega^1(x_p,x_q)$.}
\end{equation}
Then, $\J$ satisfies the Palais--Smale condition on
$\Omega^1(x_p,x_q)$.
\end{lemma}

\begin{proof}
Let $(x_k)_k \subset \Omega^1(x_p,x_q)$ be such that
\begin{equation}\label{ps2}
(\J(x_k))_k\; \mbox{is bounded}\quad \mbox{and} \quad \lim_{k \to
+\infty}\J'(x_k) = 0.
\end{equation}
As $\J$ is coercive, \eqref{ps2} implies that $(\|\dot x_k\|)_k$
is bounded; hence, there exists a compact subset $K$ in $\m_0$
such that $x_k(s) \in K$ for all $s\in I$ and all $k\in \N$.
Therefore, $(x_k)_k$ is bounded in $H^1(I,\m_0)$, thus in
$H^1(I,\R^N)$ (as $\m_0$ is isometrically embedded in $\R^N$), and
there exists $x \in H^1(I,\R^N)$ such that
\[
x_k \wk x \; \mbox{weakly in $H^1(I,\R^N)$}\;\;\mbox{and}\;\; x_k
\to x \; \mbox{uniformly in $I$}
\]
(up to subsequences). Clearly, assumption $(h_1)$ implies $x \in
\Omega^1(x_p,x_q)$.\\
Moreover, by \cite[Lemma 2.1]{BF} there exist
two sequences $(\xi_k)_{k}$ and $(\nu_k)_{k}$ in
$H^1(I,\R^N)$, with $\xi_k \in T_{x_k}\Omega^1(x_p,x_q)$, such that
\begin{equation}
\label{limits1}
\begin{split}
&x_k - x = \xi_k + \nu_k\;\;\mbox{for all $k \in \N$,}\\
&\xi_k \wk 0 \;\; \mbox{weakly and}\;\;
    \nu_k \to 0\; \;\mbox{strongly in $H^1(I,\R^N)$.}
\end{split}
\end{equation}
In order to prove that $\xi_k \to 0$ strongly in $H^1(I,\R^N)$,
consider $y_k = \phi_y(x_k)$, $t_k = \phi_t(x_k)$ and $z_k =
(x_k,y_k,t_k)$. As the coefficients $A$, $B$, $C$ in
\eqref{metric} are bounded in $K$, and $\acca$ in \eqref{stima1}
is bounded far away from zero in $K$, then so are the sequences
$(a(x_k))_k$, $(b(x_k))_k$ and $(c(x_k))_k$. Whence, from
\eqref{uai}, \eqref{ti} and \eqref{ps3}, it follows that also
$(\dot y_k)_k$ and $(\dot t_k)_k$ are  bounded in $L^2(I,\R)$.
From \eqref{newfunct} and \eqref{ps2} it follows
\[
\J'(x_k)[\xi_k]\ =\ f'(z_k)[(\xi_k,0,0)] = o(1)  ,
\]
i.e.,
\[
\begin{split}
o(1)\ =\ & \int_0^1 \langle\dot x_k,\dot \xi_k\rangle\; {\rm d}s
+\ \frac{1}{2}\ \int_0^1 \langle \nabla A(x_k),\xi_k\rangle\ \dot y_k^2\; {\rm d}s \\
&+ \int_0^1 \langle \nabla B(x_k),\xi_k\rangle\ \dot y_k \dot t_k\; {\rm d}s
- \frac{1}{2}\ \int_0^1 \langle \nabla C(x_k),\xi_k\rangle\ \dot t_k^2\; {\rm d}s .
\end{split}
\]
So, (\ref{limits1}) and previous remarks give
\[
\begin{split}
& \int_0^1 \langle\nabla A(x_k),\xi_k\rangle\ \dot y_k^2 \; {\rm
d}s\ =\ o(1) ,\quad
\int_0^1 \langle\nabla B(x_k),\xi_k\rangle\ \dot y_k \dot t_k \; {\rm d}s\ =\ o(1) ,\\
&  \int_0^1 \langle\nabla C(x_k),\xi_k\rangle\ \dot t_k^2 \; {\rm d}s\ =\ o(1),\quad
\int_0^1 \langle \dot x,\dot \xi_k\rangle\; {\rm d}s = o(1) ,\quad
\int_0^1 \langle\dot \nu_k,\dot \xi_k\rangle\; {\rm d}s = o(1).
\end{split}
\]
In conclusion, we obtain $\displaystyle \int_0^1 \langle\dot
\xi_k,\dot \xi_k\rangle {\rm d}s = o(1)$, which completes the
proof.
\end{proof}

Now, we are ready to give the proofs of our main results.

\begin{proof}[Proof of Theorem \ref{connect1}]
Fix $z_p = (x_p,y_p,t_p)$, $z_q = (x_q,y_q,t_q) \in \m$, with
$x_p$, $x_q \in \m_0$ and $(y_p,t_p)$, $(y_q,t_q) \in \R^2$. By
hypotheses $(h'_2)$ and $(h'_3)$, the case $(i)$ in the table
above holds for all $x \in \Omega^1(x_p,x_q)$ (recall
\eqref{stima1}, \eqref{definition1}). Moreover,
\begin{equation}\label{g}
\lambda_+(x)\ \geq\ \frac{a(x)-c(x)}2\ >\ 0.
\end{equation}
Thus, by (\ref{newact}), (\ref{g}), (\ref{definition1}), the
expression of $\gamma$ in $(h'_3)$ and \eqref{differe}, it is
\begin{eqnarray*}
\J(x) &\geq& \frac 12\ \|\dot x\|^2 \ -\
\frac{\Delta_-^2(x)}{a(x)-c(x)}
\\
&=&
\frac 12\ \|\dot x\|^2
\ -\
\Delta_-^2(x)
\ \left(\int_0^1 \frac {A(x)-C(x)}{\acca(x)}\; {\rm d}s\right)^{-1}
\\
&=&
\frac 12\ \|\dot x\|^2
\ -\
\Delta_-^2(x)
\ \left(\int_0^1 \frac {1}{\gamma(x)}\; {\rm d}s\right)^{-1}
\\
&\geq& \frac 12\ \|\dot x\|^2 \ -\ \sqrt{\Delta_y^2+\Delta_t^2} \
\left(\int_0^1 \frac {1}{\gamma(x)}\; {\rm d}s\right)^{-1}.
\end{eqnarray*}
Define
\[
\bar J(x)\ :=\ \frac 12\ \|\dot x\|^2 \ -\ \sqrt{\Delta_y^2+\Delta_t^2} \
\left(\int_0^1 \frac {1}{\gamma(x)}\; {\rm d}s\right)^{-1}.
\]
Note that $\bar J$ has the same form of the static functional $J$
in \eqref{newfunc}. Moreover, from \eqref{stima1}, the scalar
field $\gamma$ is strictly positive and satisfies hypothesis
\eqref{ipotesi1}, which is analogous to condition $(H_2)$ in
Theorem \ref{maintheorem}. So, from Proposition \ref{inf}, it
follows that $\J$ is bounded from below and coercive. Furthermore,
Lemma \ref{ps1} implies that $\J$ satisfies the Palais--Smale
condition on $\Omega^1(x_p,x_q)$. Thus, Theorem \ref{min} applies,
and a geodesic connecting $z_p$ with $z_q$ exists.
\end{proof}

\begin{proof}[Proof of Theorem \ref{connect2}]
Fix $z_p = (x_p,y_p,t_p)$, $z_q = (x_q,y_q,t_q) \in \m$, with
$x_p$, $x_q \in \m_0$ and $(y_p,t_p)$, $(y_q,t_q) \in \R^2$. By
hypotheses $(h''_2)$ and $(h''_3)$, the case $(iv)$ in the table
above holds for all $x \in \Omega^1(x_p,x_q)$. Then, functional
$\J$ in (\ref{newact}) is not only bounded from below, but also
coercive, as
\[
\J(x) \ \ge\ \frac 12\ \|\dot x\|^2\qquad \hbox{for all $x \in \Omega^1(x_p,x_q)$.}
\]
Then, by hypothesis $(h_1)$, Lemma \ref{ps1} applies, and $\J$
satisfies the Palais--Smale condition on $\Omega^1(x_p,x_q)$.
Finally, as $\Omega^1(x_p,x_q)$ is also complete, Theorem
\ref{min} implies the existence of a critical point for $\J$ on
$\Omega^1(x_p,x_q)$. Hence, a geodesic connecting $z_p$ with $z_q$
exists.
\end{proof}

The same arguments which allow one to prove the global property of
the geodesic connectedness as stated in Theorems \ref{connect1}
and \ref{connect2}, can be used for proving the existence of a
geodesic joining two fixed points.

\begin{proposition}
Let $(\m= \m_0 \times \R^2,\langle\cdot,\cdot\rangle_L)$ be a
G\"odel type spacetime such that $(h_1)$ holds and fix two points
$z_p = (x_p,y_p,t_p)$, $z_q = (x_q,y_q,t_q) \in \m$ such that
${\cal L}(x)\geq\nu>0$ (resp. ${\cal L}(x)\leq -\nu<0$) for all
$x\in \Omega^1(x_p,x_q)$. If $(h'_3)$ (resp. $(h''_3)$) holds,
then $z_p$ and $z_q$ are geodesically connected.
\end{proposition}

\begin{remark}
{\rm If case $(ii)$ occurs for all $x \in \Omega^1(x_p,x_q)$, the
opposite inequality for the difference $a(x)-c(x)$ prevents
to proceed as in the proof of Theorem \ref{connect1} (cf. \eqref{g}).\\
On the other hand, if case $(iii)$ occurs for all $x \in \Omega^1(x_p,x_q)$, then
$\lambda_+(x) \ge \lambda_-(x) > 0$. Then, from \eqref{newact}
\[
\J(x)\geq \frac 12\ \|\dot x\|^2 - \frac 12 \frac{(\Delta_+^2(x) +
\Delta_-^2(x))}{\lambda_-(x)}.
\]
Clearly, it is possible to give suitable conditions for
$\lambda_-(x)$ on $\Omega^1(x_p,x_q)$ which ensure the
coercivity of $\J$.
Nevertheless, the expression of $\lambda_-(x)$ makes hard the
analytic formulation of these bounds.}
\end{remark}

\bigskip

\centerline{\bf Acknowledgment}

\smallskip

\noindent
J.L. Flores would like to acknowledge the Dipartimento di Matematica, Universit\`a degli
Studi di Bari ``Aldo Moro'', where this work started, for its kind
hospitality.


\end{document}